\documentclass[11pt,a4paper]{article}

\usepackage{latexsym}
\usepackage{amsfonts}

\newcommand{\noin}{\noindent}

\newtheorem{theorem}{\bf Theorem}[section]
\newtheorem{lemma}[theorem]{\bf Lemma}
\newtheorem{corollary}[theorem]{\bf Corollary}

\newtheorem{remark}[theorem]{\bf Remark}
\newtheorem{definition}[theorem]{\bf Definition}

\begin{document}

\title{Eigenvalue estimates for submanifolds \\ with locally  bounded mean  curvature  }
\author{G. Pacelli Bessa \and J. F\'{a}bio Montenegro  }
\date{\today}
\maketitle
\begin{abstract}We give lower bounds estimates for the first Dirichilet eigenvalues for domains $\Omega $ in  submanifolds $M\subset N$ with locally bounded mean curvature. These lower bounds depend   on the local injectivity radius, local upper bound for sectional curvature  of $N$  and local bound for the mean cuvature of $M$. For sumanifolds with bounded mean curvature of Hadamard manifolds these lower bounds depends only on the dimension and the bound on mean curvature.
\end{abstract}
\section{Introduction}Let $M$ be a   Riemannian manifold and  $\Omega \subset M$ a connected domain  with compact closure and nonempty boundary $\partial \, \Omega $. The first eigenvalue of $\Omega$ denoted by $\lambda_{1}(\Omega )$ is given by
\begin{equation}\label{eq1}
\lambda_{1}(\Omega )=\inf\left\{\frac{\int_{\Omega}\vert \nabla f\vert^{2}}{\int_{\Omega} f^{2}},\,f\in {L^{2}_{1,0}(\Omega )\setminus\{0\}}\right\},
\end{equation} where $L^{2}_{1,0}(\Omega )$ is the completion  of $C^{\infty}_{0}(\Omega )$ with respect to the norm
\begin{equation}
\vert\varphi \Vert^2=\int_{\Omega}\varphi^{2}+\int_{\Omega} \vert\nabla \varphi\vert^2.\label{eq2}
\end{equation}   Observe that in this variational formulation, the boundary $\partial \,\Omega $ needs not to be smooth.  It is clear that if $\Omega_{1}\subset \Omega_{2}$ are bounded domains, then $ \lambda_{1}(\Omega_{1})\geq \lambda_{1}(\Omega_{2})\geq 0.$ From this considerations one may define the first eigenvalue of a complete non-compact Riemannian manifold $M$ as follows. 
 \begin{definition}The first eigenvalue of a complete non-compact Riemannian manifold $M$ is 
 \begin{equation}\lambda_{1}(M)=\lim_{r\rightarrow \infty}\lambda_{1}(B_{M}(p,r)),\label{eq3}\end{equation} where $B_{M}(p,r)$ is the geodesic ball of radius $r$ centered at $p$.
\end{definition}
This limit (\ref{eq3})  exists and it is independent on the center $p$. 
 Since $\lambda_{1}(M)$ can be zero, it is natural to ask   under what conditions a complete non-compact Riemannian manifold $M$ has positive first eigenvalue. Toward this question,  Mckean  proved the following theorem.
\begin{theorem}[Mckean, \cite{kn:mckean}] Let $M$ be an $m$-dimensional, complete non-compact, simply connected Riemannian manifold with sectional curvature $K_{M}\leq -a^2<0$, then
\begin{equation}\lambda_{1}(M)\geq \frac{(m-1)^{2}a^{2}}{4}\label{eq4}
\end{equation}
\end{theorem}

Leung-Fu Cheung and Pui-Fai Leung in \cite{kn:c-l} considered non-compact $m$-dimensional submanifolds $M$ of the hyperbolic space $\mathbb{H}^{n}$ with bounded mean curvature $H$ and showed that if $\vert H \vert \leq c $ for some constant $c<m-1$, then $\lambda_{1}(M)\geq (m-1-c)^{2}/4$. They considered the distance  function $r$ (distance to  a point) in the hyperbolic space  $\mathbb{H}^{n}(-1)$ and showed  that when $r$ is restricted to a non-compact submanifold $M$  with mean curvature $\vert H\vert \leq c <n-1$ then $\Delta\, r\geq n-1-c$, where $\Delta $ is the laplace operator in the induced metric on $M$.  After that, using cleverly Young inequality concluded that $\lambda_{1}(M)\geq (n-1-c )^{2}/4.$ 

 In this paper we will consider similar questions regarding eigenvalues estimates. Let $M\subset_{\varphi} N$  be an immersed submanifold with locally bounded mean curvature in a Riemannian manifold $N$ (see section 4 for a definition). Let $p\in N$ and consider connected components   $\Omega \subset \varphi^{-1} (B_{N}(p,r)) $. We give lower bounds estimates for the first Dirichilet eigenvalue $\lambda_{1}(\Omega)$. The radius of the ball and the lower bounds  depend on the sectional curvatures of the $N$ in $B_{N}(p,r) $, on the injectivity radius $inj(p)$ of $N$ at $p$ and the mean curvature of $\Omega $ (see Theorem (\ref{thmApp1})). When $N$ is complete, simply connected and has  sectional curvature bounded from above by a negative constant the lower bounds for  $\lambda_{1}(\Omega)\geq c>0$, where $c$ is independent on on $r$. Thus, we can conclude that $\lambda_{1}(M)>c$ (see Corollary \ref{CorApp1}). This extends Cheun-Leung result to submanifolds with bounded mean curvature of Hadamard manifolds.

\section{Main Lemma  (Lemma \protect{\ref{mainlemma}})}

 As mentioned before,  Cheung-Leung showed  that  the distance function $r$ when  restricted to a non-compact $m$-dimensional submanifold $M\subset \mathbb{H}^{n}(-1)$  with mean curvature $\vert H\vert <m-1$, the laplacian of $r$ in the induced metric on $M$ is bounded from below, $\Delta\, r\geq n-1-\sup \,\vert H \vert$. Using Young inequality they showed that $\lambda_{1}(M)\geq (n-1-\sup\,\vert H\vert )^{2}/4.$ Their idea (and proof) applies immediately to vector fields $X$ on domains $\Omega $ of manifolds provided that  $\inf \, div\, X \,>0$ and $\sup \, \vert X\vert <\infty $ on $\Omega $. This is the content of the following lemma.
\begin{lemma}Let $\Omega \subset M$ be a domain with compact closure $(\partial \Omega \neq \emptyset )$   in a Riemannian manifold $M$  and $X:M\rightarrow TM$ be a smooth nonzero vector field with
$ \sup_{\Omega } \,\vert X \vert < \infty $ and  $ \inf_{\Omega } \, div\,X>0$, then \begin{equation}\lambda_{1}(\Omega)\geq \left[\frac{\inf_{\Omega}\, div\,X}{2 \,\sup_{\Omega}\,\vert X \vert}\right]^{2}>0.\label{eqML1}\end{equation}\label{mainlemma}
\end{lemma}

\noin {\em Proof:}
Let $f\in C^{\infty}_{0}(\Omega )$. The vector field $f^{2}\,X$ has compact support on $\Omega $. Computing $div(f^{2}X)$ we have
\begin{eqnarray}div(f^{2}X) & =   & \langle grad\,f^{2}, X\rangle + f^{2}div (X)\nonumber \\
                            &\geq & -\vert grad\,f^{2}\vert\cdot \vert X \vert+ \inf \, div\,X\cdot f^{2} \label{eqML2} \\
                            &\geq & -2 \cdot \sup \,\vert X \vert
\cdot \vert f\vert\, \vert grad\,f\vert + \inf \, div\,X\cdot f^{2}\nonumber 
\end{eqnarray}
Using the inequality \begin{eqnarray}
 -2\, \vert f\vert\, \cdot \vert grad\, f\vert &  \geq & -\epsilon \,\cdot  \vert f\vert^{2} - 1/\epsilon \,\cdot \vert grad \, f \vert^{2}\nonumber \end{eqnarray} for all $\epsilon >0$, we have from (\ref{eqML2}) that 
\begin{eqnarray}
                  div(f^{2}X)          &\geq & \sup \,\vert X \vert
\,\cdot (-\epsilon\,\cdot  \vert f \vert^{2} - \frac{1}{\epsilon}\, \cdot \vert grad\,f\vert^{2})+ \inf \, div\,X \cdot f^{2} \label{eqML3}
\end{eqnarray}
Integrating (\ref{eqML3}) over $\Omega $ we have  that
\begin{eqnarray}
0=\int_{\Omega} div(f^{2}X) &\geq & \sup \,\vert X \vert \cdot
\,\int_{\Omega}(-\epsilon\, \vert f \vert^{2} - \frac{1}{\epsilon}\, \vert grad\,f\vert^{2})+ \inf \, div\,X \cdot \int_{\Omega}f^{2}, \nonumber
\end{eqnarray}
therefore
\begin{eqnarray}
  \int_{\Omega}\vert grad\,f\vert^{2}&\geq & \frac{\epsilon}{\sup \,\vert X \vert
}\,(\inf \, div\,X\,-\,\sup \,\vert X \vert\cdot \epsilon )\,\int_{\Omega} f^{2}.\nonumber
\end{eqnarray}
Choosing $\epsilon = (\inf \, div\,X)/(2\sup \,\vert X \vert)
$ we have that 
\begin{eqnarray} \int_{\Omega}\vert grad\,f\vert^{2}& \geq &  \left[\frac{\inf \, div\,X}{2\sup \,\vert X \vert
}\right]^{2}\,\,\int_{\Omega} f^{2}.\label{eqML4}
\end{eqnarray}
By the variational formulation of the first eigenvalue of $\Omega $ the inequality  (\ref{eqML4}) implies (\ref{eqML1}).
  
\begin{corollary}There are no smooth bounded vector fields $X:M\rightarrow TM$ with $\inf_{M}\,div\,X>0$ on complete non-compact manifolds with $\lambda_{1}(M)=0$. In particular there is no such vector fields on $\mathbb{R}^{n}$.
\end{corollary}

\section{Basic Formulas}  In order to aplly Lemma (\ref{mainlemma}) we need to to find suitable vector fields $X\in TN$ such that the divergent of the tangencial component $X^{t}\in TM\subset TN $ has positive infimum. The natural vector field is the gradient of the distance function that depending on the geometry of $M$ and $N$ has divergent strictly positive. To see this we shall set up some basic fomulas relating the divergent of the gradient of the distance function, its lower bound, the mean curvature and upperbounds of the sectional curvatures of $N$.

  Let $\varphi : M \hookrightarrow N$ be an isometric immersion, where $M$ and $N$ are complete Riemannian manifolds. Consider a smooth function $g:N \rightarrow \mathbb{R}$ and the composition $f=g\,\circ\,\varphi :M \rightarrow \mathbb{R}$. Identifying $X$ with $d\varphi (X)$ we have at $p\in M$ and  for every $X\in T_{p}M$ that
\begin{eqnarray} \langle grad\,f\, , \, X\rangle \,=\,d f(X)=d g(X)\,=\, \langle grad\,g\, ,\, X\rangle, \nonumber
\end{eqnarray}
therefore 
\begin{equation} grad\,g\,=\, grad\, f\, +\,(grad\,g)^{\perp},\label{eqBF1}
\end{equation}
where $(grad\,g)^{\perp}$ is perpendicular to $T_{p}M$.
Let $\nabla $ and $\overline{\nabla}$ be the Riemannian connections on $M$ and $N$ respectively,  $\alpha (p) (X,Y) $ and $Hess\,f(p)\,(X,Y)$  be respectively the second fundamental form of the immersion $\varphi $ and the Hessian of $f$ at $p\in M$,  $X,Y \in T_{p}M$. Using the Gauss equation we have that
 \begin{equation}Hess\,f (p)  \,(X,Y)= Hess\,g (\varphi (p))\,(X,Y) + \langle grad\,g\,,\,\alpha (X,Y)\rangle_{\varphi (p)}.
\label{eqBF2}
\end{equation}
 Now let $\{ e_{1},\ldots e_{m}\}$ be an orthonormal basis for $T_{p}M$.  Then  taking the trace in (\ref{eqBF2}) we have  that
\begin{eqnarray} 
\Delta \,f (p)                & = & \sum_{i=1}^{m}Hess\,f (p)  \,(e_{i},e_{i})\nonumber \\
                             & = & \sum_{i=1}^{m}Hess\,g (\varphi (p))\,(e_{i},e_{i}) + \langle grad\,g\,,\,\sum_{i=1}^{m}\alpha (e_{i},e_{i})\rangle.\label{eqBF3}
\end{eqnarray} 
 We should mention that the formulas (\ref{eqBF2}) and (\ref{eqBF3}) are  well known in the literature, see \cite{kn:c-l}, \cite{kn:j-k} , \cite{kn:c-g}. Another important tool  is the Hessian Comparison Theorem. 

\begin{theorem}[Hessian Comparison Theorem] Let $M$ be    a  complete Riemannian manifold  and  $x_{0},x_{1} \in M $.   Let
 $\gamma:[0,\,\rho (x_{1}) ]\rightarrow M$ be a minimizing geodesic joining $x_{0}$ and $x_{1}$ where $\rho (x)$ is the  distance function $dist_{M}(x_{0}, x) $. Let $K$ be the sectional curvatures of $M$  and  $\mu_{i}(\rho )$, $i=0,1$, be these functions defined below.
 \begin{equation} \mu_{0}(\rho )=\left\{ \begin{array}{lcll} & k_{0} \cdot\coth (k_{0}\cdot \rho (x)),  & if  & \inf_{\gamma} K=-k_{0}^{2} \\
                                                             &                                          &     &   \\
                                                             & \frac{1}{\rho (x)},                      & if  &  \inf_{\gamma} K=0  \\
                                                             &                                          &     & \\
                                                             & k_{0} \cdot \cot (k_{0} \cdot \rho (x)), & if  &  \inf_{\gamma} K =k_{0}^{2}\; and \; \rho < \pi/2k_{0}. 
\end{array}\right.\label{eqBF4}
\end{equation}

\vspace{.1cm}

And

\vspace{.1cm}

\begin{equation} \mu_{1}(\rho )=\left\{ \begin{array}{lcll}  & k_{1} \cdot\coth (k_{1} \cdot\rho (x)), & if  & \sup_{\gamma} K=-k_{1}^{2} \\
                                                             &                       &     &   \\
                                                             & \frac{1}{\rho (x)},   & if  &  \sup_{\gamma} K=0  \\
                                                             &                       &     & \\
                                                             & k_{1} \cdot\cot (k_{1}\cdot \rho (x)),  & if  &  \sup_{\gamma} K =k_{1}^{2}\; and \; \rho < \pi/2k_{1}. 
\end{array}\right.\label{eqBF5}
\end{equation}
\vspace{.1cm}

\noin Then the Hessian of $\rho $ and $\rho^{2}$ satisfies
\vspace{.1cm}

\begin{equation}\begin{array}{rclll}\mu_{0}(\rho(x))\cdot\vert X\vert^{2}       & \geq & Hess\,\rho(x)(X,X) & \geq & \mu_{1}(\rho(x))\cdot\vert X\vert^{2}\\
                                                                           &      &                    &      & \\
 2\rho (x)\cdot\mu_{0}(\rho(x))\cdot\vert X\vert^{2} & \geq & Hess\,\rho^{2}(x)(X,X) & \geq & 2\rho (x)\cdot\mu_{1}(\rho(x))\cdot\vert X\vert^{2}.
\end{array}\label{eqBF6}
\end{equation}
Where $X$ is any vector in $T_{x}M$ perpendicular to $\gamma'(\rho (x))$.\label{thmHess}
Taking the trace in (\ref{eqBF6}) we have 
\begin{equation} 
\begin{array}{rclll}(n-1)   \mu_{0}(\rho(x))      & \geq & \Delta \,\rho(x)& \geq & (n-1)\mu_{1}(\rho(x))\\
                                                                           &      &                    &      & \\
 2(n-1)\rho (x)\mu_{0}(\rho(x)) & \geq & \Delta\,\rho^{2}(x) & \geq & 2(n-1)\rho (x)\mu_{1}(\rho(x)).
\end{array}\label{eqBF7}
\end{equation}
Here $n= dim\, M$.

\end{theorem}

\section{Applications of the Main Lemma}

Let $N$ be a complete Riemannian manifold  and   let $\kappa (p,r)=\sup\{K_{N}(x);\, x\in B_{N}(p,r)\}$, where $K_{N}(x)$ are the sectional curvatures of $N$ at $x$ and  $B_{N}(p,r)$ is the geodesic with center $p$ and radius $r$. 
As an application of the lemma (\ref{mainlemma})  we are going to give Dirichilet eigenvalue lower bound estimate for domains  on immersed submanifolds with locally bounded mean curvature. 
\begin{definition} An isometric immersion $\varphi : M \hookrightarrow N$ has locally bounded mean curvature $H$ if  for any $p\in N$ and $r>0$ this number
$h(p,r)=\sup \{\vert H(x)\vert ;\,x\in \varphi (M)\cap B_{N}(p,r) \} $ is finite. 
\end{definition}
Cheng in (\cite{kn:cheng1}) gave lower bounds for the first Dirichilet eigenvalue on balls $B_{N}(p,r)$, $r<inj (p)$. These bounds are the first Dirichilet eigenvalue on balls of the same radius on space forms with constant curvature $\sup_{B_{N}(p,r)} K_{N}$. 
The following theorem can be viewed as a  version of Cheng's  eigenvalue estimates for balls in submanifolds. 
 \begin{theorem}\label{thmApp1}
Let 
$\varphi : M \hookrightarrow N$ be an isometric immersion with locally bounded mean curvature $ H$ where $M$ is an  $m$-dimensional complete non-compact Riemannian manifold. 
Choose any point $p\in N\setminus \varphi (M)$ and $r$ as follows:

\noin Case I. The injectivity radius of $N$ at $p$ is finite, $inj(p)<\infty$.

\begin{enumerate}

\item If  $ \kappa (p,inj(p))=k^2,$ choose  $r<\,\min\,\{ inj(p),\,\pi/2k,\, \cot^{-1} [\,h(p,inj(p))/(m-1)k]/k\,\}$. 

\item If $  \kappa (p,inj(p))=0$, choose $\,\,r<\,\min\,\{ inj(p), (m-1)/h(p,inj(p))\}$.  

\item If $\kappa (p,inj(p))=-k^{2}$ and $h(p,inj(p))<m-1$, then  choose $r<inj (p)$. 

If $h(p,inj (p))\geq (m-1)$ then choose $r<\,min\,\{inj(p)\,,\,\coth^{-1} [\,h(p,inj(p))/(m-1)k]/k\,\}$

\end{enumerate}

\noin Case II.  The injectivity radius of $N$ at $p$ is infinity, $inj(p)=\infty$.

\noin \hspace{.32cm} 4.\hspace{.02cm} If $\lim_{r\rightarrow \infty}\kappa (p,r)=\infty$,
choose $r=\max_{s>0}\min\,\{ \frac{\pi}{2\sqrt \kappa (p,s)}\,,\, \cot^{-1}[\,\frac{h(p,s)}{(m-1)k(p,2)}]/k(p,s)\}$.

\noin \hspace{.32cm} 5.\hspace{.02cm} If $\lim_{r\rightarrow \infty}\kappa (p,r)=0,-k^{2}$, choose $r$ as in 2. and 3. respectively.

\vspace{.3cm}

\noin Let $\Omega$ be  any connected component of 
$ \varphi^{-1}(\, \overline{B_{N}(p,r)}\,)$, $r$ chosen above.
Then the first eigenvalue $\lambda_{1}(\Omega )$  for the Dirichilet problem on $\Omega $  has  the following lower estimates.

\begin{itemize}

\item $r$ as in  1. and 4.: $\lambda_{1}(\Omega )\geq [(m-1)\,k\,\cot \, (k\,r)-h(p,r)]^{2}/4$, where in 4. we substitute $k$ by $\sqrt (k(p,r)$.

\item $r$ as in  2.: $\lambda_{1}(\Omega )\geq [(m-1)/r \, -h(p,r)]^{2}/4 $

 \item $r$ as in 3.: If $h(p,inj(p))<m-1$ then $\lambda_{1}(\Omega )\geq [(m-1)-h(p,r)]^{2}/4 $. 

\noin If $h(p,inj(p))\geq m-1$ then $\lambda_{1}(\Omega )\geq [(m-1)k\coth (k\, r)-h(p,r)]^{2}/4 $.

\end{itemize}
\end{theorem}

\noin {\em Proof:} Let $\rho (x)=dist_{N}(p,\,x\,)$ be the distance  function on $N$ and $f=\rho\circ \varphi :M\rightarrow \mathbb{R}$.   $f$ is  smooth on $\varphi ^{-1}(B_{N}(p,inj(p)))$. Consider a connected component $\Omega $ of $ \varphi^{-1}(\, \overline{B_{N}(p,r)}\,)$ and in this conponnent let $X=\,grad\,f$. Now by (\ref{eqBF3}), $div\,X (x)   = \sum_{i=1}^{m-1} \,Hess \, \rho\,(\varphi (x))(e_{i},e_{i}) + \langle grad \,\rho \,,\, H )\rangle (\varphi (x))$, $\{ e_{1},\ldots , e_{m}\}$ base orthonormal of $T_{x}M$. Applying the Hessian Comparison Theorem we have that
\begin{itemize}
\item $div \,X \geq (m-1)\,k\,\cot \, (k\,r)-h(p,r)>0$, if $ \kappa (p,r)=k^2$.
\item $div \,X \geq  (m-1)/r \, -h(p,r)>0$, if $  \kappa (p,r)=0$.
\item $div \,X \geq  m-1)\,k\,\coth \, (k\,r)-h(p,r)>0$, if $\kappa (p,r)=-k^{2}$.
\end{itemize}
Since $\vert X\vert \leq 1$ the estimatives follows from Lemma (\ref{mainlemma}).

\begin{remark} It was proved by Jorge and Xavier in \cite{kn:j-x} (Lemma 1) that if $\varphi : M\rightarrow N$ is an isometric immersion with locally bounded second fundamental form $\alpha $ and $r$ chosen as in Theorem (\ref{thmApp1}) where $h(p,r)$ is substituted by $\beta (p,r)=\{\vert \alpha (x)\vert ;\,x\in \varphi (M)\cap B_{N}(p,r) \} $, then
$ \varphi^{-1}(\, \overline{B_{N}(p,r)}\,)$ is a countable collection  of disjoint  compact topological $m$-disks $\Omega_{j} $ and the function $f=\rho\circ \varphi $ has only one critical point $q_{j}\in \Omega_{j}$ and its a minimum.
Moreover, $diam_{M}(\Omega_{j})$ is bounded above by $C(r,\, \rho (\varphi (q_{j})))$. It also shwon that $\Omega_{j}$ contains a geodesic ball centered at $q_{j}$ with radius $r_{j}\geq c(r,\rho (\varphi (q_{j})), \inf_{B_{N}(p,r)}\, K_{N})$. 

\end{remark}

\begin{corollary}Let $\varphi : M \hookrightarrow N$ be an isometric immersion with  bounded mean curvature $H$, where $M$ is    an $m$-dimensional complete non-compact Riemannian manifold $M$    and $N$ is an $n$-dimensional complete  simply connected Riemannian manifold  $N$ is with sectional curvature $K_{N}$ satisfying

 $$K_{N}\, \leq \,-\,a^{2} \, < \,0, \label{eqAppB}\;\; {\rm for\; a\, positive \; constant } \; \, a \, >\,0.$$

If  $\vert H \vert\leq \beta < (m-1)\,a $, then 
\begin{equation}\lambda_{1}(M)\geq \frac{[(m-1)a-\beta]^{2}}{4}>0\label{eqAppC}.\end{equation} In particular, there is entire Green's functions on $M$. \label{CorApp1}
\end{corollary}

\noin {\em Proof:} As in Theorem (\ref{thmApp1}),
 consider $\rho=dist_{N}(p,\,\cdot\,)$, $\;p\in N\setminus \varphi (M)$,  $f=\rho\circ \varphi $ and $X=grad\,f$. Since $N$ is simply connected and is negatively curved, every point of $N$ is a pole, hence $inj (p)=\infty$. The estimative  for the case $inj(p)=\infty $ and $h(p,r)<m-1$ for all $r>0$ in Theorem (\ref{thmApp1}) does depend on $r$, i.e. $\lambda_{1}(\Omega )\geq (m-1 -\beta)^{2}/4$.
 Since 
 any domain $\Omega \subset M$ is immersed inside a ball $B_{N}(p,r)$ for sufficient large $r$ we have that $\lim \lambda_{1}(\Omega_{j})\geq  (m-1 -\beta)^{2}/4$ for any exhaustion $\{\Omega_{j}\subset \Omega_{j+1}\}$ of $M$.
 The existence of entire Green's functions on $M$ is  a diretc consequence of Theorem A.3 from \cite{kn:s-y}, page 84.

\end{document}